\documentstyle[12pt,epsfrag,psfrag]{amsart}

\begin{document}

\title[Retardation of Plateau-Rayleigh Instability]{Retardation of Plateau-Rayleigh Instability:\\
        a Distinguishing Characteristic\\Among Perfectly Wetting Fluids}
\author{John McCuan}
\address{\hskip-\parindent
McCuan, Mathematical Sciences Research Institute, 1000 Centennial Drive, 
Berkeley, CA  94720}
\email{mccuan@@msri.org}
\thanks{This work was supported in 
 part by NASA grant NAS8-39225 (Gravity Probe B Relativity Mission).  
Research at MSRI is supported in part by NSF grant DMS-9022140.}

\begin{abstract}
We consider a cylindrical film of fluid adhering to a rigid cylinder of 
fixed radius.  The main result is to give the critical (maximum) length for 
which such a film of given thickness can be stable.  The problem is 
considered both when the cylinder remains stationary and when the fluid and 
the cylinder are co-rotating at a fixed angular velocity.  
\end{abstract}

\maketitle

\section{Young's Equation}

\begin{figure}[ht]
\centerline{{\psfrag{x}{$x$}
                   \psfrag{g}{$\gamma$}
                   \psfrag{t}{$\gamma$}
                   \psfrag{p1}{$\Pi_1$}
                   \psfrag{p2}{$\Pi_2$}
         \epsfysize=2.5in
         \epsfxsize=3.5in
         \leavevmode\epsfbox{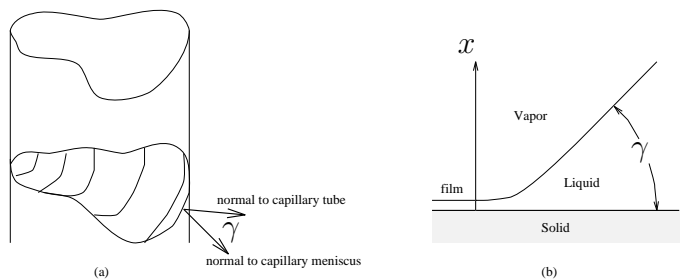}}}
\centerline{ }
\caption{Contact Angles\label{p-angles}}
\end{figure}

The equation of Young asserts that the angle at which a fluid surface in 
equilibrium meets homogeneous walls containing the fluid is constant 
(Figure~\ref{p-angles}(a)). 
This assertion has been the subject of debate 
since it was introduced; Adamson \cite[pg. 360]{AdaPhy} writes 
concerning the microscopic contour of a meniscus,
``The likely picture for the case of a nonwetting drop on a flat surface 
is shown in Figure~\ref{p-angles}(b).  There is a region of negative 
curvature as the drop profile joins the plane of the solid.''  
Presumably, such a film (and perhaps the region of negative curvature) 
are not subject to variation.  Nevertheless, if such a region exists even  
microscopically, then a mathematical description of 
Figure~\ref{p-angles}(b) incorporating the macroscopic contact angle is 
of obvious interest.  We offer such an explanation in an (admittedly) 
degenerate 
case for which the film {\em is} subject to variation and (for 
this reason) does play a role on the macroscopic level.  

\section{Perfectly Wetting Fluids}

Films of certain fluids assert their presence conspicuously due to their 
ability to transport the fluid.  We have, here, helium at low 
temperature in mind, but there are other examples.  From the point of 
view of equilibria, however, these fluids have typically been treated 
under the assumptions that the film is negligible and that Young's 
condition prevails with $\gamma = 0$.  Considering the typical thicknesses 
of such films (20-90 atomic diameters for helium) these assumptions 
may seem justified.  A recent experiment \cite{NesDem} involving rotating 
nearly cylindrical films indicates otherwise.  We now give a brief 
description of this experiment.

\begin{figure}[ht]
\centerline{{\psfrag{a}{$\alpha$}
                   \psfrag{g1}{$\gamma_1$}
                   \psfrag{g2}{$\gamma_2$}
                   \psfrag{p1}{$\Pi_1$}
                   \psfrag{p2}{$\Pi_2$}
         \epsfysize=1in
         \epsfxsize=2in
         \leavevmode\epsfbox{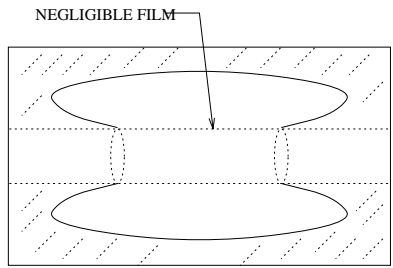}}}
\centerline{ }
\caption{Axially Symmetric Doubly Connected Bubble\label{spin}
\label{spin-can1}}
\end{figure}

In low gravity \footnote{For this discussion we may assume zero 
gravity; in the actual experiment low gravity was simulated on the 
earth's surface.} an annular container consisting of the region between 
two coaxial cylinders was partially filled with a perfectly wetting 
fluid.  In this instance the inner cylinder was made of steel and the fluid 
was silicone oil.  Previous applications of the Young condition 
$\gamma = 0$ and the assumptions described above predicted that 
under rotation about the axis of the cylinders (at certain speeds 
with certain volumes of fluid) a stable axially symmetric equilibrium 
would be observed with cross-section resembling that shown in 
Figure~\ref{spin-can1}.

What was in fact observed was no equilibrium at all, but the nearly 
cylindrical film on the inner cylinder wall thickened periodically and 
formed ``pendant drops'' which were flung outward, returning to the 
bulk of fluid near the outer cylinder.  In this way the fluid was 
``pumped'' by the capillary instability of the ``film'' on the inner 
cylinder in a cycle apparently against the centripetal force in part of 
its course.  See Figure 3.
\begin{figure}[ht]
\centerline{{\psfrag{a}{$\alpha$}
                   \psfrag{g1}{$\gamma_1$}
                   \psfrag{g2}{$\gamma_2$}
                   \psfrag{p1}{$\Pi_1$}
                   \psfrag{p2}{$\Pi_2$}
         \epsfysize=2.5in
         \epsfxsize=3.5in
         \leavevmode\epsfbox{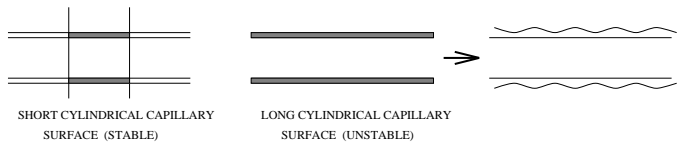}}}
\centerline{ }
\caption{Thickening of Cylindrical Film\label{thicken}}
\end{figure}

It is this experiment which prompted the work to be presented below.

\section{A New Approach}

In the standard Young-Laplace theory of equilibrium capillary surfaces, 
for which the reader is referred to \cite{FinEqu}, equilibrium surfaces 
$\cal{S}$ and their stability are determined by an energy functional of the 
form
\begin{equation}\label{star}
	E=\sigma|{\cal S}|-\sigma\beta|{\cal W}|+\int_X{\cal U}
\end{equation}
where $\sigma$ is the surface tension, $|{\cal S}|$ is the area of the 
free surface, $|{\cal W}|$ is the area of the wetted region, $X$ is the 
volume of fluid, and ${\cal U}$ is a potential depending on position; 
the Young equation is given by $\cos \gamma = \beta$ where $\beta$ is 
the ``adhesion coefficient.''  One might attempt to 
distinguish between various ``zero contact angle'' fluids in this 
framework by choosing an adhesion coefficient $\beta > 1$.  It has 
recently been shown by B. White \cite{WhiExi} that such an approach 
yields no difference (in equilibria nor stability) from the case 
$\beta = 1$.

In the perfectly wetting case we suggest that the second term in 
(\ref{star}) be dropped altogether and a ``Van der Waals'' potential 
be incorporated in ${\cal U}$.  For containers with suitable geometry 
this additional potential depends only on the distance $d$ to the 
container and, in any event, is required to become infinite and negative 
as $d \searrow 0$.  In this way, the issue of a contact line or a contact 
angle need never be addressed because all equilibria must wet every 
portion of the container.  

In the remainder of the paper we carry out an analysis of certain 
equilibria in a particularly simple container geometry.  The results 
are rather surprising particularly in view of the historical introduction of 
the next section.

\section{Plateau-Rayleigh Instability}

Plateau in 1873 determined experimentally that the length at which a 
cylindrical column of liquid becomes unstable is between 3.13 and 3.18 
times the diameter of the column.\footnote{This experiment was repeated in
\cite{MasExp} by G.~Mason who obtained a value $3.143 \pm .004$.}  
The theoretical value of $\pi$ for this multiple has been given by 
Rayleigh and more recently by Barbosa and do~Carmo \cite{RayIns,BarSta}.  
We present a simple (and partial) derivation based on the second 
variation of (\ref{star}).  In this 
case $|{\cal W}|={\cal U}\equiv 0$, and the first variation $\delta E$ 
gives rise to the equation of constant mean curvature 
$2H=\lambda = {\rm constant}$ to which the cylinder is a solution with 
$H=1/2r$.  

Calculating the second variation about the cylindrical solution one 
obtains
\begin{equation}
	\delta^2 E=\left(-{1\over{r^2}}\phi_{\theta\theta}-\phi_{zz}-
			{1\over{r^2}}\phi \right)\phi
\end{equation}
where $\theta$ and $z$ are cylindrical coordinates on the cylinder 
${\cal S}$.  Thus, using the criterion of Barbosa and DoCarmo
\cite{BarSta}, the cylinder will be unstable when $L \phi = -(1/r^2)
\phi_{\theta\theta}-\phi_{zz}-(1/r^2)\phi$ has a non-positive 
eigenvalue whose eigenfunction satisfies 
\begin{equation}\label{vol}
	\int_{\cal S} \phi=0.
\end{equation}
This latter condition is the  result of constrained volume.  Solving 
$L \phi=\mu \phi$ by a separation of variables $\phi=A(\theta)B(z)$ 
with the boundary conditions $A(0)=A(2\pi),\ A'(0)=A'(2\pi);\ 
B(0)=0=B(l)$ one obtains eigenvalues
\begin{equation}
	\mu_{k,m}={k^2 -1\over{r^2}}+{m^2 \phi^2\over{l^2}}
\end{equation}
for $k=0,1,2,\ldots$ and $m=1,2,3,\ldots$.

$\mu_{k,m} >0$ for $k>0$.  Thus, instability can only arise for 
eigenvalues $\mu_{0,m}$.  The smallest of these is 
$\mu_{0,1}=\pi^2/l^2 -1/r^2$, but the corresponding eigenfunction 
$\phi_{0,1}=\sin (\pi/l)z$ does not satisfy (\ref{vol}) and must be 
discarded.  The next smallest is $\mu_{0,2}=4\pi^2/l^2 - 1/r^2$ 
whose eigenfunction does satisfy (\ref{vol}) and gives rise to an 
instability exactly when 
\begin{equation}\label{threestars}
	l \ge 2\pi r.
\end{equation}

This is a symmetric variation, and at the length indicated in 
(\ref{threestars}) the column 
of fluid will develop a neck which becomes smaller in radius until 
the column separates into two pieces.  It is of note that in the 
experiment of G.~Mason \cite{MasExp}, the column of fluid was 
centered on a silica rod and 
the necking apparently proceeded until the column ``ruptured'' around the 
rod.  According to our model such behavior is excluded for perfectly 
wetting fluids.  

\section{Rotating Cylindrical Films}

We now carry over the analysis of the previous section to determine the 
critical length for the onset of instability of a cylindrical film 
of perfectly wetting fluid on a rotating, coaxial, rigid cylinder.

\begin{figure}[ht]
\centerline{{\psfrag{r}{$r_0$}
                   \psfrag{s}{$r_0$}
                   \psfrag{g}{$\omega>0$}
                   \psfrag{z}{$\omega=0$}
         \epsfysize=2.5in
         \epsfxsize=3.5in
         \leavevmode\epsfbox{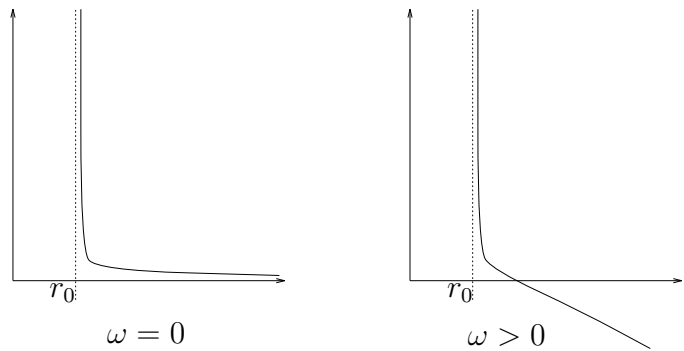}}}
\centerline{ }
\caption{The Shape of $\partial{\cal U}/\partial{r}$\label{Uprime}}
\end{figure}

The rotational potential energy is given by ${\cal U}_1=-(\rho\omega^2/2)
r^2$, and we take a model Van der Waals potential approximating 
that for helium ${\cal U}_2=-\alpha/d^3$ where $d=r-r_0$.  
A general form for 
$L \phi$ has been calculated by Wente and others \cite{WenThe}:
\begin{equation}
	L \phi=-\Delta \phi -2(2H^2-K)\phi+\nabla{\cal U}\cdot N \phi
\end{equation}
where $\Delta$ is the intrinsic Laplacian on the free surface ${\cal S}$, 
$H$ and $K$ are the mean and Gauss curvatures of ${\cal S}$ and $N$ is 
the outward normal to ${\cal S}$.  From this one obtains by separation 
of variables as before
\begin{equation}
	\mu_{k,m}={\partial{\cal U}\over{\partial r}}+{k^2 -1\over{r^2}}+
		{m^2 \pi^2\over{l^2}}
\end{equation}
for $k=0,1,2,\ldots$ and $m=1,2,3,\ldots$.  The two smallest contributing 
eigenvalues are 
\begin{equation}
	\mu_{0,2}={\partial{\cal U}\over{\partial{r}}}-{1\over{r^2}}+
			{4\pi^2\over{l^2}}
\end{equation}
corresponding to the axisymmetric variation of the last section and 
\begin{equation}
	\mu_{1,1}={\partial{\cal U}\over{\partial{r}}}+{\pi^2\over{l^2}}
\end{equation}
corresponding to a non-axisymmetric variation.  The typical shape of 
$\partial{\cal U}/\partial{r}$ is shown in Figure~\ref{Uprime} and evidently 
results in three cases.

\medskip

\noindent VERY THIN FILMS:  There is a value $r_1$ determined by the equation
\begin{equation}\label{cross}
	{\partial{\cal U}\over{\partial r}}(r_1)={1\over{r_1^2}}
\end{equation}
such that any cylindrical film (of arbitrary length) satisfying 
$r_0<r<r_1$ is stable.

For helium in the absence of rotation (\ref{cross}) becomes 
\begin{equation}
	{3\alpha\over {(r-r_0)^4}}={1\over{r^2}}.
\end{equation}
Thus, $r=r_{\rm{critical}}$ satisfies
\begin{equation}
	r_{\rm{critical}}=r_0 +{\sqrt{3\alpha}\over{2}}+ {1\over{2}}
		\sqrt{4r_0 \sqrt{3\alpha}+3\alpha}.
\end{equation}

\medskip

\noindent MEDIUM THICKNESS FILMS:  Without rotation all cylindrical films with 
$r_1 \le r$ are unstable to the axially symmetric variation associated 
to $\mu_{0,2}$ whenever 
\[
	l \ge l_{0,2}\equiv {2\pi\over{\sqrt{{1\over{r^2}}-
		{\partial{\cal U}\over{\partial r}}}}}.
\]

In the presence of rotation there is a value $r_2>r_1$ satisfying the 
equation 
\begin{equation}
	{\partial {\cal U}\over{\partial r}}(r_2)=-{1\over{3r_2^2}}.
\end{equation}
If $r_1 \le r < r_2$, then an {\em axially symmetric\/} instability sets in 
at $l=l_{0,2}$.

\medskip 

\noindent THICK FILMS (only in the presence of rotation):  If $r_2\le r$, 
then the cylinder of radius $r$ becomes unstable to the 
{\em non-axisymmetric\/} variation associated to $\mu_{1,1}$ when $l$ 
exceeds
\[
	l_{1,1}={\pi\over{\sqrt{-{\partial{\cal U}\over{\partial{r}}}}}}.
\]
NOTE:  The only properties of the virtual (Van der Waals) potential 
${\cal U}_2$ which have been used are
\begin{description}
\item{(0)}  $\lim_{r\searrow r_0}{\cal U}_2(r)=+\infty$.

\item{(i)}  \[  {\partial^2{\cal U}_2\over{\partial{r}^2}}<0.  \]

\item{(ii)}  (When $\omega =0$)  ${\partial{\cal U}\over{\partial{r}}}$ 
decays faster than $1/r^2$ as $r\to\infty$.
\end{description}

\bibliographystyle{plain}
\bibliography{/chern0/john/papers/bib/books,/chern0/john/papers/bib/papers,/chern0/john/papers/bib/heresay}

\end{document}